# Home health care planning with considering flexible starting/ending points and service features

**Pouria Khodabandeh [a], Vahid Kayvanfar [a], Majid Rafiee [a], Frank Werner [b],[*]**

[a] *Department of Industrial Engineering, Sharif University of Technology, Tehran, Iran*

[b] *Faculty of Mathematics, Otto-von-Guericke-University, Magdeburg, Germany*

**ABSTRACT**

One of the recently proposed strategies in health systems is providing services to patients at home, improving the service quality, besides reducing the health system costs. In the real world, some services, such as biological tests or blood sampling, force the nurses to start or end his/her route from/at the laboratory instead of the depot, changing the whole optimal planning. The effect of these special service requirements and features has not been considered so far. In this study, a new mathematical model is suggested considering the flexibility of starting/ending places of each nurse's route according to the specific characteristics of each service. Then several sets of problems in various sizes are solved using the proposed model, where the results confirm the efficiency of the proposed approach. In addition, some sensitivity analyses are performed on the parameters of the required features of the services, followed by some managerial insights and directions for future studies.

**Keywords:** Home health care; Integrated routing and scheduling; Optimization; Flexible starting and ending points; Special service requirements.

## 1 INTRODUCTION

Health has always been one of the most important human concerns from the past, and there has been a lot of efforts to create stable health in communities. Today, with an increase in life expectancy and a decline in birth rates, in most societies, we see a rise in the average age of the population, which has imposed a lot of costs on health systems around the world, and a significant share of countries' budget is spent on health care.

In the current situation, creating a plan to use the resources optimally seems to be necessary due to resource limitations. For example, human resource shortages and shortage of hospitals are a limiting factor in health services. On the other hand, one can also see an individual life pattern of elderly people in the developed and developing countries that makes it necessary to focus on this group of individuals. Developing a framework for offering health care to individuals at their home effectively could be a good solution to this problem. In this context, home health care (HHC) companies are usually faced with conflicts in their targets such as maximizing the offering service level and minimizing the operating costs which usually should be simultaneously optimized. With the attendance of various companies in this field, new angles of existing costs have been clarified, which have encouraged researchers to use optimization techniques in this field. One of the most related issues for companies consists in finding an optimal plan for offering services. This problem is called the home health care routing and scheduling problem (HHCRSP) which has attracted a lot of attention recently among the researchers.

Fernandez, Gregory, Hindle, and Lee (1974) were the first to study routing and scheduling in the HHC context. They investigated a working day of community nurses in their research.

---

[*] Corresponding author.
E-mail address: frank.werner@mathematik.uni-magdeburg.de

Hindle, Hindle, and Spollen (2000) and Hindle, Hindle, and Spollen (2009) extended the problem of Fernandez et al. (1974) by considering resource allocations and travel cost estimation aspects. Decision support systems (DSSs) in this context were first considered by Bertels and Fahle (2006) and Eveborn, Flisberg, and Rönnqvist (2006). This problem was developed as an extension of the vehicle routing problem by Akjiratikarl, Yenradee, and Drake (2007) and solved by using a particle swarm optimization (PSO) algorithm.

In recent years, the HHCRSP has attracted the attention of many researchers, and various studies have been implemented in this field. Recent research in this field can be classified into two categories: deterministic and uncertain studies. In the category of recent deterministic researches, Mankowska, Meisel, and Bierwirth (2014) provided a daily planning with regard to the specific needs of each patient, the particular skills of each nurse, as well as the relationships between different services. Liu et al. (2013) applied the research of Mankowska et al. (2014) to the logistics domain of home health care. They noted that nurses may need to transfer medicines and medical supplies from the company center to the patients' home and return the laboratory samples, unused drugs, and medical equipment to the center. R. Liu, Xie, and Garaix (2014) extended the problem of the multi-period VRP with time windows (PVRPTW) to three different patient demand types.

Decerle, Grunder, El Hassani, and Barakat (2018) addressed the problem of home health care with considering a few objectives simultaneously and focusing on the applicability of the planning. They employed a memetic algorithm to solve this problem. Fathollahi-Fard, Hajiaghaei-Keshteli, and Tavakkoli-Moghaddam (2018) presented a bi-objective green home health care problem that considers environmental pollution. The conventional HHCRSP is extended to demand and capacity management by J. Nasir and Dang (2018). A mixed integer programming (MIP) model is suggested with considering workload balancing. To solve this model, a variable neighborhood search (VNS) algorithm as well as a heuristic method were applied. J. A. Nasir, Hussain, and Dang (2018) developed a mathematical model to consider group-based and telehealth-based care services. Three different goals were considered in their model to achieve this target. The first goal was to select an optimal location for the HHC centers, staff and group patients' centers considering required specifications. The second goal of their study was to schedule the group patients' sessions in a specific time windows by pairing each nurse with a telephone service staff. Finally, as their third goal, the violation from good quality services and the dissatisfaction of the patients were considered.

Fikar and Hirsch (2018) introduced the concept of car and trip sharing within home health care nurses. They also considered the option that the nurses are allowed to walk to the patients' home. Lin, Hung, Liu, and Tsai (2018) presented a modified harmony search (MHS) algorithm that considers three problems of nurse rostering, nurse routing and nurse re-rostering simultaneously. They specified that the past works had considered the issues of the nurse rostering problem (NRP) and VRPTW independently, while it has been investigated at the same time in this study.

The home health care problem in Chinese communities was considered by Zhang, Yang, Chen, Bai, and Chen (2018). These communities have an intense distribution of patients. M. Liu, Yang, Su, and Xu (2018) developed a bi-objective model to minimize the company costs and to improve patient satisfaction. In their study, multiple weeks of planning and medical teams were considered. To solve such a problem, the epsilon-constraint method was used to obtain Pareto fronts of the problem.

Shanejat-Bushehri, Tavakkoli-Moghaddam, Momen, Gasemkhani, A., and Tavakkoli-Moghaddam (2019) considered the HHCRSP problem with temporarily precedence and synchronization constraints as well as limited allowable times for transferring the collected biological samples to the laboratory. The goal was to minimize the cost related to the transportation and the idle time of the caregivers. They presented a mathematical model and applied simulated annealing and tabu search in two phases. Grenouilleau, Legrain, Lahrichi, and Rousseau (2019) presented a method for the HHCRSP problem which is based on a set partitioning formulation as well as a variable neighborhood search framework. Their algorithm solved first a linear relaxation, and then a constructive heuristic was applied to generate an integer solution. Euchi, Zidi, and Laouamer (2019) presented a distributed optimization approach for the HHCRSP problem, which used artificial intelligence techniques. They integrated automatic learning and search techniques to optimize the assignment of caregivers to patients. Kohdabandeh, Kayvanfar, Rafiee and Werner (2021) included downgrading aspects into the classical HHCRSP problem. They considered the additional goal to minimize also the difference between the actual and the potential skills of the nurses. To solve the derived bi-objective model, an epsilon-constraint approach was applied. They also performed a sensitivity analysis on the epsilon parameter. Ghiasi, Yazdani, Vahdani, and Kazemi (2021) considered the HHCRSP problem with two transportation modes, namely public and private modes under the multi-depot version. The objective was the minimization of the sum of the travel distance and overtime costs. After presenting a mixed integer programming model, three metaheuristic algorithms have been given for solving large instances. Xiang, Li, and Szeot (2021) considered the HHCRSP problem with minimizing total costs and maximizing patients' preference satisfaction. They formulated a bi-objective mixed integer linear programming model. For solving this problem, a local search algorithm was embedded into the basic framework of a nondominated sorting genetic algorithm. The algorithm obtained approximate Pareto-optimal solutions for small instances in a shorter computation time than the epsilon-constraint method.

In the category of recent uncertain studies, Yuan, Liu, and Jiang (2015) addressed the HHCRSP with stochastic service times. The problem was modeled in the form of a stochastic programming problem with recourse, in which the expected value of late arrivals of the nurses was considered. Then the branch and price (B&P) method was used to solve this problem. In R. Liu, Yuan, and Jiang (2018), stochastic service and travel times were considered in their problem. They employed a chance constraint in order to guarantee the probability of on-time services. In their study, a route-based mathematical model was developed, and the branch and cut (B&C) algorithm was used with the discrete estimation method to solve it. Moreover, labeling algorithms and acceleration methods were also employed for solving the proposed model. Lanzarone and Matta (2014) presented a robust strategy for home health care, in which random patient requests were considered, and the nurse allocation to the patient was investigated as well. Rodriguez, Garaix, Xie, and Augusto (2015) modeled the HHC problem considering a stochastic demand for patients and handled it through stochastic programming. Then the B&C method was used to solve their proposed model. Shi, Boudouh, and Grunder (2017) addressed the HHC problem with a fuzzy demand for the patients. They employed a fuzzy chance constraint method as well as a hybrid genetic algorithm (GA) for the solution of the model. Shi, Boudouh, Grunder, and Wang (2018) investigated the HHC problem considering stochastic travel and service times through stochastic programming with recourse. Then they applied a simulated annealing (SA) algorithm to get the solutions of the problem. Issabakhsh, Hosseini-Motlagh, Pishvaee, and Saghafi Nia (2018) presented a robust optimization model for patients in need of dialysis, in which the travel times of nurses are uncertain. An innovative approach to accepting or rejecting new patients by a nurse was recently investigated in Demirbilek, Branke, and Strauss (2019). Random scenarios were created to present this simple and fast method. Carello, Lanzarone, and Mattia (2018) proposed a multi-criteria optimization approach to consider different goals of the stakeholders with regard to the need for continuity of the care and considering certain and uncertain patient demands. Khodaparasti, Bruni, Beraldi, Maleki, and Jahedi (2018) presented a multi-period allocation method through a robust approach in which the patient demand was considered to be uncertain.

Bazirha, Kadrani, and Benmansour (2020) dealt with the HHCRSP problem with stochastic travel and care times with the goal to minimize the transportation costs of the caregivers and the expected value of recourse caused by delayed services and the overtime of the caregivers. The performance of the developed genetic algorithm with an embedded Monte Carlo simulation is discussed. Later in an extended paper, Bazirha, Kadrani, and Benmansour (2021) proposed a stochastic programming model with recourse for the HHCRSP problem, where the goal was to minimize the transportation cost and the expected value of recourse and also multiple services and their synchronization were considered. While the underlying deterministic problem was solved by CPLEX, a genetic algorithm and a general variable neighborhood search heuristic, the stochastic problem was solved by Monte Carlo simulation embedded into their genetic algorithm. Recently, Bazirha (2022) applied a similar approach to the HHCRSP problem with additional hard/fixed time windows and the goal to minimize the traveling costs of the caregivers and the average number of unvisited patients. Very recently, Di Mascolo, Martinze and Espinouse (2021) presented a literature survey and a bibliometric analysis in the field of home health care. They reviewed and analyzed the current state-of-the-art with a focus on uncertain and dynamic aspects.

To the best of the authors' knowledge and as it can be seen from the discussed literature, most of the researches assumed that each nurse should start her/his journey from the depot and end it at the laboratory. In home health care's real world, there are some services which need special instruments that forces the nurse to start her/his route from the laboratory. On the other hand, if a nurse has to take biological test or a blood test she/he should end her/his route at the laboratory instead of depot. Such these requirements led us to develop conventional models to a flexible model which can consider different options for origin and destination of the nurses' routes.

The rest of the paper is organized as follows. The description of the problem and the mathematical model are given in Section 2. Section 3 discusses the computational experiments. The results of a sensitivity analysis are explained in Section 4. Section 5 offers some managerial insights and finally, some conclusions and future studies are presented in Section 6.

## 2 PROBLEM DESCRIPTION

One of the most important issues that has always been a focus of attention in the HHCRSP is to make the problem conditions closer to the real world of the health industry. In home health care, there are many different types of services with different features that play an important role in the real world. Therefore, they should always be considered by the planners. Each of these services needs to take into account their own particular considerations according to their characteristics. On the other hand, each of the nurses has special skills and to provide the appropriate services to any patient, these services should be offered with a set of these skills.

In home health care, many patients require several different services, which may be provided by a nurse or several nurses with relevant skills. One of the most important assumptions existing in the home health literature was that each patient needs only one service. Such a problem with this assumption is very simple. However, it does not seem to be valid in the real world since some patients need several services. The services provided to patients in the field of home health care are very diverse, including services with specific characteristics that are required to satisfy their anterior and posterior requirements. Some services require a specific equipment, materials, and supplies that the nurse must obtain at the start of her/his journey from the lab, while other services do not require such an action. Moreover, some services require an equipment that must be delivered to the laboratory upon the completion of the nurse's route, or some vital services, such as biological tests, are required to be delivered to the laboratory at the end of the day. Due to such real-world requirements, it seems to be necessary to create a flexibility on the starting and ending points of the nurse's route. This flexibility is determined depending on the specific features of the services.

In this study, a mathematical model is proposed with the consideration of flexible starting and ending points of the routes that are determined depending on the features of the required services. This problem has also been developed by a simultaneous consideration of different services for the patients, the patient time windows, the correct sequence of the patient services as well as the rest of the other constraints required for the VRP problems. The effectiveness of this proposed novel model in the optimal planning and scheduling process is demonstrated through solving the samples problems.

*2.1 Assumptions*

- Each patient may need several services, where all of these needs should be satisfied by the company through its nurses.
- The starting and ending positions of each nurse's route (depot or laboratory) are determined by the features of the services that are included in the route.
- Every patient has a desirable time window that must be respected.
- Each patient's service could begin after the end of the previous patient's service in addition to the time that is needed to travel to the place of the second patient.
- The travel time, the service time and the demand of each patient is deterministic and the planner knows them before the beginning of the planning.
- The planning is performed daily.
- Each nurse has a car for traveling between the patients, and multi-mode traveling and travel sharing concepts are not considered.
- Urgent service calls and emergent situations are not considered in the planning.

*2.2 Notations*

*2.2.1 Subscripts*

| | |
|---|---|
| $i$ | Starting point of each transfer ($i=1,2,..., n+1$), where $n$ is the number of patients. |
| $j$ | Ending point of each transfer ($j=2,3,...n+2$); where $n$ is again the number of patients. |
| $k$ | Nurse index ($k = 1,2,…,V$), where $V$ is the number of nurses. |
| $s$ | Service index ($s = 1,2,…,S$), where $S$ is the number of different services in the planning. |

*2.2.2 Sets*

| | |
|---|---|
| $C$ | Set of patients. |
| $N$ | Set of all nodes that includes the depot, the patients and the laboratory. |
| $V$ | Set of nurses. |
| $S$ | Set of services. |

*2.2.3 Input Parameters*

| | | |
|---|---|---|
| $t_{ij}$ | Traveling time from node *i* to node *j*. | |
| $t_{is}$ | Needed time for giving service *s* to patient *i*. | |
| $l_i$ | Minimum acceptable value of the time window of patient *i*. | |
| $u_i$ | Maximum acceptable value of the time window of patient *i*. | |
| $a_{ks}$ | Nurse qualification matrix; 1 means that nurse *k* has the qualification of service *s*. | |
| $g_{js}$ | Patient's service needs matrix; 1 means that service *s* is needed by patient *j*. | |
| $R_{s1}$ | Service starting requirements matrix; 1 means that service *s* needs to start from the laboratory. | |
| $R_{s2}$ | Service ending requirements matrix; 1 means that service *s* needs to end at the laboratory. | |

### 2.3 Decision Variables

| | |
|---|---|
| $x_{ijks}$ | 1 if nurse *k* travels from node *i* to node *j* for giving service *s*; 0 otherwise. |
| $S_{iks}$ | Starting time of giving service *s* to patient *i* by nurse *k*. |
| $\delta_{k1}$ | 1 if nurse *k* needs to start her/his route from the laboratory; 0 otherwise. |
| $\delta_{k2}$ | 1 if nurse *k* needs to finish her/his route at the laboratory; 0 otherwise. |

### 2.4 The mathematical model

#### 2.4.1. Objective function

$$Min \; z = \sum_{i \in N} \sum_{j \in N} \sum_{k \in V} \sum_{s \in S} t_{ij} \cdot x_{ijks} \tag{1}$$

Constraints

$$\varepsilon \cdot \sum_{i \in N} \sum_{j \in N} \sum_{s \in S} R_{s1} x_{ijks} \leq \delta_{k1} \leq \sum_{i \in N} \sum_{j \in N} \sum_{s \in S} R_{s1} x_{ijks} \quad \forall k \in V \tag{2}$$

$$\sum_{i \in C} \sum_{s \in S} x_{1iks} + M \cdot \delta_{k1} \geq 1 \quad \forall k \in V \tag{3}$$

$$\sum_{i \in C} \sum_{s \in S} x_{1iks} - M \cdot \delta_{k1} \leq 1 \quad \forall k \in V \tag{4}$$

$$\sum_{i \in C} \sum_{s \in S} x_{(n+2)iks} + M \cdot (1 - \delta_{k1}) \geq 1 \quad \forall k \in V \tag{5}$$

$$\sum_{i \in C} \sum_{s \in S} x_{(n+2)iks} - M \cdot (1 - \delta_{k1}) \leq 1 \quad \forall k \in V \tag{6}$$

$$\varepsilon \cdot \sum_{i \in N} \sum_{j \in N} \sum_{s \in S} R_{s2} x_{ijks} \leq \delta_{k2} \leq \sum_{i \in N} \sum_{j \in N} \sum_{s \in S} R_{s2} x_{ijks} \quad \forall k \in V \tag{7}$$

$$\sum_{i \in C} \sum_{s \in S} x_{i1ks} + M \cdot \delta_{k2} \geq 1 \quad \forall k \in V \tag{8}$$

$$\sum_{i \in C} \sum_{s \in S} x_{i1ks} - M \cdot \delta_{k2} \leq 1 \quad \forall k \in V \tag{9}$$

$$\sum_{i \in C} \sum_{s \in S} x_{i(n+2)ks} + M \cdot (1 - \delta_{k2}) \geq 1 \quad \forall k \in V \tag{10}$$

$$\sum_{i \in C} \sum_{s \in S} x_{i(n+2)ks} - M \cdot (1 - \delta_{k2}) \leq 1 \quad \forall k \in V \tag{11}$$

$$\sum_{i \in N} \sum_{s \in S} x_{ijks} - \sum_{i \in N} \sum_{s \in S} x_{jiks} = 0 \quad \forall j \in C, k \in V \tag{12}$$

$$S_{iks_1} + t_{is_1} + t_{ij} - M(1 - x_{ijks_2}) \leq S_{jks_2} \quad \forall i,j \in N, k \in V, s_1 \in S, s_2 \in S \tag{13}$$

$$l_i \leq S_{iks} \leq u_i \qquad \forall i \in C, k \in V, s \in S \qquad (14)$$

$$\sum_{k \in V} \sum_{i \in N} a_{ks} x_{ijks} = g_{js} \qquad \forall j \in C, s \in S \qquad (15)$$

$$x_{ijks} = a_{ks} \cdot g_{js} \qquad \forall i, j \in N, k \in V, s \in S \qquad (16)$$

Equation (1) presents the optimization criterion of the model, i.e., to minimize the total traveling time that the nurses of the company spend for providing suitable care to the patients. Constraints (2) to (6) are used to determine the starting place of the routes of the nurses. The initial position of each nurse is decided with respect to the features of all services that should be applied in that route. If there is a service in a route which forces the caregiver to start her/his journey from the laboratory instead of the depot, the variable $\delta_{k1}$ will become 1 by constraint (2). In this situation, constraints (3) and (4) are deactivated, and constraints (5) and (6) are activated and guarantee that the nurse must start her/his route from the laboratory. Contrariwise, constraints (3) and (4) are activated and enforce the nurse to start her/his journey from the depot. Constraints (7) to (11) are used to determine the ending place of the routes of the nurses. The finishing position of each nurse is decided with respect to the features of all services that should be applied in that route. If there is a service in a given route, which forces the caregiver to finish her/his journey from the laboratory instead of the depot, the variable $\delta_{k2}$ will become 1 by constraint (7). In this situation, constraints (8) and (9) are deactivated, and constraints (10) and (11) are activated and guarantee that the nurse must finish her/his route at the laboratory. Inversely, constraints (8) and (9) are activated and enforce the nurse to finish the journey at the depot. Constraint (12) ensures that each nurse should depart from the patient's place after giving care and go to another patient's home. Constraint (13) states that a service could be initiated after ending the previous service as well as the needed time for traveling the nurse to the new place. Constraint (14) indicates that each patient's service should be started in her/his desirable time window. Constraint (15) is used to ensure that, if patient $j$ needs service $s$, exactly one of the nurses with the required qualifications should go to the patient's home and serve her/him. Constraint (16) ensures that for giving service $s$ by nurse $k$ to patient $j$, the patient must need the service $s$ and of course, the nurse must have the required qualifications.

## 3 COMPUTATIONAL EXPERIMENTS

In this study, in order to demonstrate the effectiveness of the proposed model in the real world, first a small example is implemented and the results are then explained accompanied by a discussion regarding the verification of the proposed model. Next, some benchmark instances that are obtained from (Mankowska et al., 2014) are applied to the model, and IBM ILOG CPLEX Optimization Studio Version 12.6.0.0 is used to solve the model. In this study, a computer with an Intel i7-4710HQ processor, 8 GB of RAM and 2.5 GHz core speed is used for the computational experiments.

### 3.1 Model verification

In this subsection, first a small instance is solved to demonstrate the effectiveness of the proposed model and to clarify the benefits of this model in comparison to the traditional home health care routing and scheduling models. The characteristics of the solved problem are shown in Table 1.

There are 3 nurses and 10 patients in this problem. The nurses have different qualifications and on the other hand, the patients have different service needs. There are 6 different services and qualifications, where service #3 needs finishing the route at the laboratory and service #6 needs starting the route from the laboratory.

*Table 1. Characteristics of the considered example*

| Characteristics of the example | | |
|---|---|---|
| *Number of patients* | *Number of nurses* | *Number of services* |
| 10 | 3 | 6 |

In this section, to show the importance of the new model compared to the old home health care routing and scheduling models, the problem without considering flexibility in the starting and ending points is first solved, and then it is solved using the proposed model in this paper. Finally, the differences and benefits of the proposed novel model are presented. The results of the classic and the proposed model are given in Table 2 and Table 3, respectively. For a better illustration of the results of the new model, the optimal routes for each nurse are demonstrated in Figure 1.

*Table 2. Optimal routes without considering flexibility in the model*

| | Optimal routes of the nurses without considering flexible starting and ending points |
|---|---|
| Nurse1 | Depot →(S6) 9 →(S1) 10 → Depot |
| Nurse2 | Depot →(S3) 11 →(S5) 3 →(S3) 6 →(S3) 8 →(S4) 5 → Depot |
| Nurse3 | Depot →(S5) 9 →(S6) 11 →(S5) 7 →(S2) 4 →(S4) 2 →(S3) 10 → Depot |

*Table 3. Optimal routes considering flexibility in the model (the proposed model)*

| | Optimal routes of the nurses considering flexible starting and ending points |
|---|---|
| Nurse1 | Lab →(S6) 9 →(S1) 10 → Depot |
| Nurse2 | Depot →(S5) 9 →(S3) 11 →(S5) 3 →(S3) 6 →(S4) 10 →(S3) 8 → Lab |
| Nurse3 | Lab →(S5) 7 →(S6) 11 →(S2) 4 →(S4) 2 →(S4) 5 → Depot |

By comparing Table 2 and Table 3, one can obviously observe that in addition to the existence of starting and ending points from and to the laboratory, the whole routes of the nurses are changed. Nurse #1 started her/his route from the depot and after giving service to patients 9 and 10 her/his route ended at the laboratory. In fact, because of the existence of service 6, the route was started from the laboratory. The initial place of nurse #2 was the depot and her/his route ended at the laboratory after giving service to patients 9,11,3,6,10 and 8. Nurse #3 started her/his route from the laboratory and after giving service to patients 7, 11, 4, 2 and 5 her/his route ended at the depot.

As it can be observed from these new results, nurse #2 visited patients 9 and 10 instead of patient 5 in the old planning. This new and different planning is owing to attendance of service 3 in the patient's journey. In addition, one can see that the route of nurse #3 was changed and she/he visited patient 5 instead of patients 9 and 10.

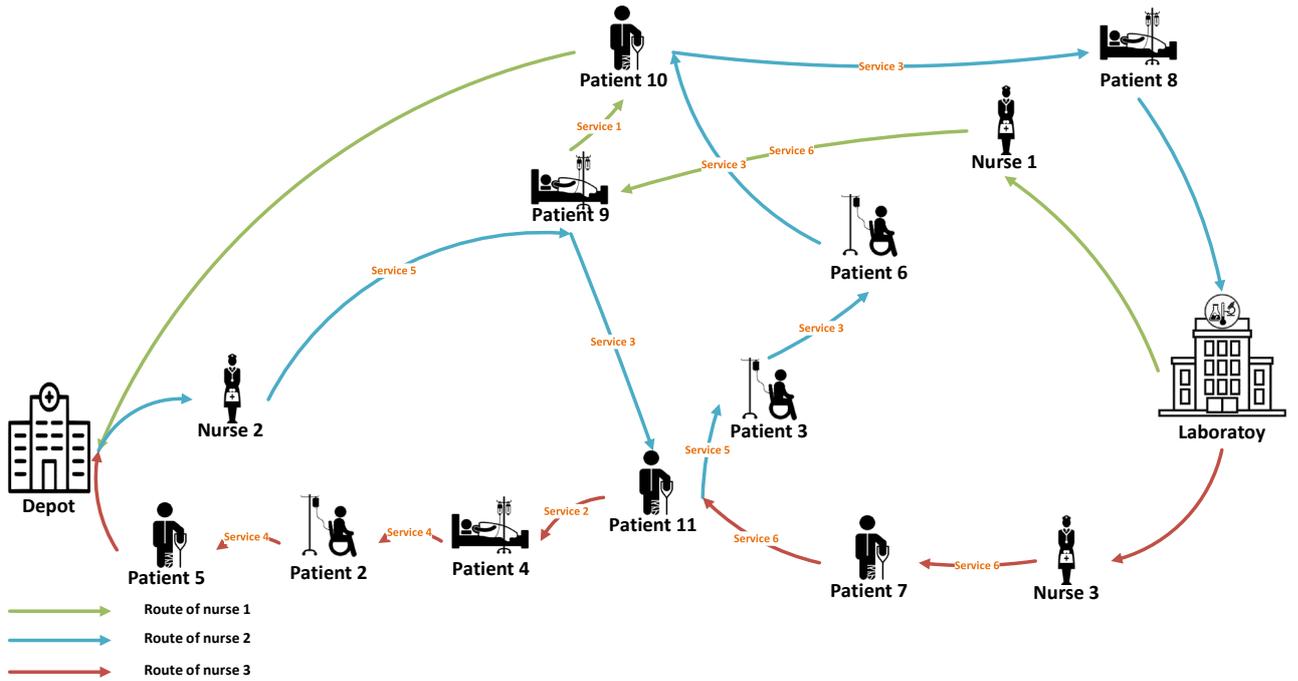

*Figure 1. Optimal routes of the nurses by applying the flexible model*

Considering the differences between the flexible proposed model and the traditional model based on the obtained results of this example, one can conclude that the service features that force our planning to start from or end at different places could be managed and taken into account by applying the novel model. These different service features affect not only individual services having these features, but also strongly affect the whole journey of each nurse. In fact, some nurses should visit different patients compared to the old planning.

### *3.2 Results*

In this subsection, three different categories of instances are employed to demonstrate the effectiveness of the suggested model. In Table 4, the characteristics of the benchmark instances are presented.

*Table 4. Characteristics of the benchmark instances*

| **Characteristics of the instances** | | | | |
|---|---|---|---|---|
| *Instance type* | *Instance number* | *Number of patients* | *Number of nurses* | *Number of services* |
| Small | #1 − #5 | 10 | 3 | 6 |
| Medium | #6 − #10 | 15 | 5 | 6 |
| Large | #11 − #15 | 25 | 5 | 6 |

To the best of the authors' knowledge, this study is the first to consider different features of services in the routing and scheduling home health care. Accordingly, the value of parameter *R* is assumed by the authors. The service requirements of starting from the depot or laboratory and finishing at the depot or laboratory are summarized in Table 5.

The obtained solutions of solving the small-sized, medium-sized and large-sized instances are presented in Table 6, Table 7 and Table 8, respectively.

*Table 5. Starting and ending features of the services*

| | **Starting and ending requirements of the services** | |
|---|---|---|
| *Services* | *Needs to start from lab* | *Needs to finish at lab* |
| *Service 1* | | |
| *Service 2* | | |
| *Service 3* | | * |
| *Service 4* | | |
| *Service 5* | | |
| *Service 6* | * | |

*Table 6. The results for the small instances*

| *Instance number* | *Optimal function value* | *Nurses with starting point of depot* | *Nurses with starting point of laboratory* | *Nurses with ending point of depot* | *Nurses with ending point of laboratory* | *Computational time(seconds)* |
|---|---|---|---|---|---|---|
| #1 | 654.596 | Nurse1 | Nurse2, Nurse3 | Nurse2, Nurse3 | Nurse1 | 0.43 |
| #2 | 724.654 | Nurse1, Nurse2 | Nurse3 | Nurse2, Nurse3 | Nurse1 | 0.28 |
| #3 | 533.423 | Nurse1, Nurse2 | Nurse3 | Nurse2, Nurse3 | Nurse1 | 2.06 |
| #4 | 568.63 | Nurse1 | Nurse2, Nurse3 | Nurse2, Nurse3 | Nurse1 | 2.56 |
| #5 | 596.976 | Nurse1 | Nurse2, Nurse3 | Nurse2, Nurse3 | Nurse1 | 1.91 |

*Table 7. The results for the medium instances*

| Instance number | Optimal function value | Nurses with starting point at the depot | Nurses with starting point at the laboratory | Nurses with ending point at the depot | Nurses with ending point at the laboratory | Computational time (seconds) |
|---|---|---|---|---|---|---|
| #6 | 1079.436 | Nurse1, Nurse2, Nurse4 | Nurse3, Nurse5 | Nurse1, Nurse3, Nurse4, Nurse5 | Nurse2 | 3.73 |
| #7 | 623.756 | Nurse1, Nurse2 | Nurse3, Nurse4, Nurse5 | Nurse1, Nurse3, Nurse4, Nurse5 | Nurse2 | 6.55 |
| #8 | 623.756 | Nurse1, Nurse2 | Nurse3, Nurse4, Nurse5 | Nurse1, Nurse3, Nurse4, Nurse5 | Nurse2 | 6.51 |
| #9 | 880.249 | Nurse1, Nurse2, Nurse3, Nurse5 | Nurse4 | Nurse2, Nurse3, Nurse4, Nurse5 | Nurse1 | 2.65 |
| #10 | 964.402 | Nurse1, Nurse2 | Nurse3, Nurse4, Nurse5 | Nurse1, Nurse3, Nurse4, Nurse5 | Nurse2 | 4.61 |

*Table 8. The results for the large instances*

| Instance number | Optimal function value | Nurses with starting point at the depot | Nurses with starting point at the laboratory | Nurses with ending point at the depot | Nurses with ending point at the laboratory | Computational time (seconds) |
|---|---|---|---|---|---|---|
| #11 | 1244.058 | Nurse1, Nurse2 | Nurse3, Nurse4, Nurse5 | Nurse1, Nurse3, Nurse4, Nurse5 | Nurse2 | 22 |
| #12 | 1291.087 | Nurse1, Nurse2 | Nurse3, Nurse4, Nurse5 | Nurse2, Nurse3, Nurse4, Nurse5 | Nurse1 | 9.92 |
| #13 | 1058.956 | Nurse1, Nurse2 | Nurse3, Nurse4, Nurse5 | Nurse1, Nurse3, Nurse4, Nurse5 | Nurse2 | 19.01 |
| #14 | 977.501 | Nurse1, Nurse2, Nurse4 | Nurse3, Nurse5 | Nurse1, Nurse3, Nurse4, Nurse5 | Nurse2 | 15.73 |
| #15 | 978.898 | Nurse1, Nurse2, Nurse5 | Nurse3, Nurse4 | Nurse1, Nurse3, Nurse4 | Nurse2, Nurse5 | 2880 |

It can be understood from the mentioned results that this model can be well applied for different daily planning of home health care organizations in different sizes and can assist them in their decisions as well.

# 4 SENSITIVITY ANALYSIS

In this section, a sensitivity analysis on the special parameter of the studied problem is carried out in order to yield a better vision about the effects of changing an important parameter over the outputs. The parameter R is chosen for the sensitivity analysis, which is directly related to the flexibility of the planning in the starting and ending points of routes. In each experiment, it is assumed that the parameter R has a specified value for service #4 and then the results for the routes of nurse #2 are explained. Different optimal routes for each experiment are shown in Table 9.

*Table 9. Sensitivity analysis on parameter R for small instances*

| Experiment number | Number of patients | Feature of service #4 | Optimal function value | Route of nurse #2 |
|---|---|---|---|---|
| #1 | 10 | | 561.813 | D → 11 (S3) → 3 (S5) → 6 (S3) → 8 (S3) → 5 (S4) → D |
| #2 | 10 | Needs to start from lab | 536.422 | L → 7 (S5) → 11 (S3) → 3 (S5) → 6 (S3) → 8 (S3) → 5 (S4) → D |
| #3 | 10 | Needs to end at lab | 609.901 | D → 9 (S5) → 11 (S3) → 3 (S5) → 6 (S3) → 5 (S4) → 8 (S3) → L |

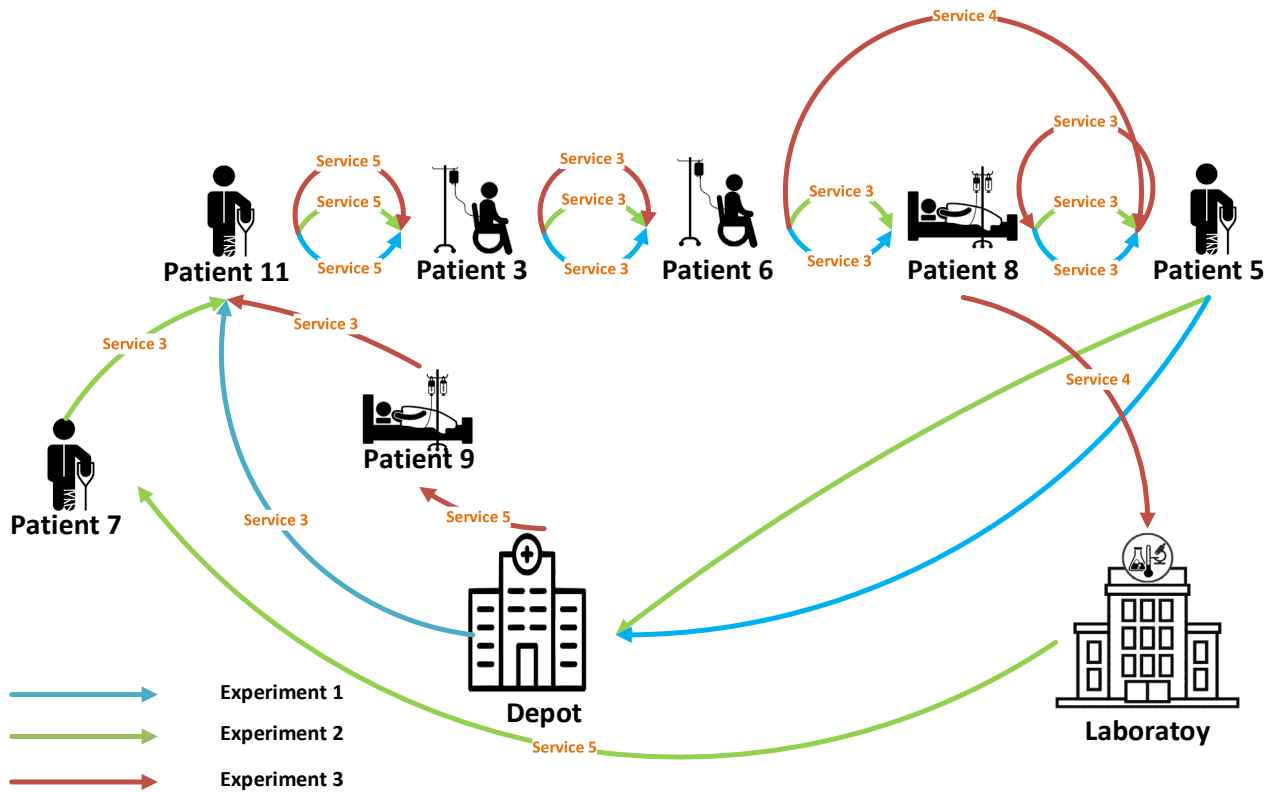

*Figure 2. Different results of the experiments for the sensitivity analysis on parameter R*

The optimal routes obtained for different experiments are also described in Figure 2.

As it is obvious from Figure 2, the whole journey of the nurse is affected through changing the starting and ending points of routes by modifying the features of service #4. In the second experiment, where service #4 forces nurse #2 to start her/his path from the laboratory, giving a service to patient #7 is also included because of changing the starting point of the nurse. In addition to experiment #2, in experiment #3, owing to the attendance of the ending special feature of service #4, servicing to patient #9 is included in the route. The obtained results of the conducted sensitivity analysis on parameter R for medium and large-sized instances are shown in Table 10 and Table 11, respectively.

According to the obtained results of the conducted sensitivity analysis on parameter R for medium-sized instances, it can be inferred that the starting and ending needs of service #4 affect only the routes of nurses #4 and #5, which is the consequence of the existing service #4 in their planned routes.

*Table 10. Sensitivity analysis on parameter R for medium instances*

| Experiment Number | Number of patients | Feature of service #4 | Optimal function value | Nurses with starting point of depot | Nurses with starting point of laboratory | Nurses with ending point of depot | Nurses with ending point of laboratory | Computational time (second) |
|---|---|---|---|---|---|---|---|---|
| **#1** | **15** | | 1050.69 | Nurse1, Nurse2, Nurse3, Nurse4, Nurse5 | | Nurse1, Nurse2, Nurse3, Nurse4, Nurse5 | | 4.09 |
| **#2** | **15** | Needs to start from lab | 1139.019 | Nurse1, Nurse2, Nurse3 | Nurse4, Nurse5 | Nurse1, Nurse2, Nurse3, Nurse4, Nurse5 | | 4.00 |
| **#3** | **15** | Needs to end at lab | 925.09 | Nurse1, Nurse2, Nurse3, Nurse4, Nurse5 | | Nurse1, Nurse2, Nurse3 | Nurse4, Nurse5 | 3.46 |

*Table 11. Sensitivity analysis on parameter R for large instances*

Sensitivity analysis on parameter R for large instances

| Experiment Number | Number of patients | Feature of service #4 | Optimal function value | Nurses with starting point of depot | Nurses with starting point of laboratory | Nurses with ending point of depot | Nurses with ending point of laboratory | Computational time (second) |
|---|---|---|---|---|---|---|---|---|
| #1 | 25 | | 1139.338 | Nurse1, Nurse2, Nurse3, Nurse4, Nurse5 | | Nurse1, Nurse2, Nurse3, Nurse4, Nurse5 | | 7.41 |
| #2 | 25 | Needs to start from lab | 1051.981 | Nurse2, Nurse3 | Nurse1, Nurse4, Nurse5 | Nurse1, Nurse2, Nurse3, Nurse4, Nurse5 | | 7.28 |
| #3 | 25 | Needs to end at lab | 1117.272 | Nurse1, Nurse2, Nurse3, Nurse4, Nurse5 | | Nurse2, Nurse3 | Nurse1, Nurse4, Nurse5 | 7.23 |

For the large instances, the special features of service #4 affect the planned routes of nurses #1, #4 and #5.

## 5 MANAGERIAL INSIGHTS

Some important managerial insights can be obtained from this study as follows:

1) In most of the traditional HHCRSP models, the planners assume that each patient requires only one service. If a patient requires three services, these service needs are considered as three different patients who have the same health profile and home location. In this study, every patient is a unique entity that can have a few service needs. This point of view can prevent patients' data redundancy. Decision makers can use this clean data as a valuable resource to make their complex decisions more efficient and manage their service and employee capacities.

2) Decision makers of home health care systems can use this novel model to have a better planning for their nurses considering different service features that may force a caregiver to have different starting and ending points. This consideration can omit changes in the planning at the day of the execution. These probable changes may happen, since some special service features may be ignored at the moment of planning by the planners in this area.

## 6 CONCLUSIONS AND FUTURE STUDIES

Nowadays, an effective management of health care systems is one of the most important concerns of policy makers. These systems can constrain high costs to the communities and have an impressive impact on the public health of the societies. At present, the correct management of health system capacities has received much attention. Giving service to the patients at their homes is one of the recent methods to provide suitable health care services. The most important issue that is addressed at the operational level of this problem by the researchers is the routing and scheduling of home health care. In order to improve the quality of the services and reduce the operational costs, it is important to find an optimal solution for the planning (including routing and scheduling) of the home health care problem. This problem is actually an extension of the famous vehicle routing problem (VRP) considering additional health care required features. This study presents a mathematical model to consider flexibility in the starting and ending positions of the caregivers due to different real-world service features. These real-world service features can affect the whole optimal planning of the nurses. In this regard, various constraints of home health care problems were taken into account such as patients' time windows, an appropriate sequence of services, the necessity of responding to all patient needs and matching the nurses' skills with the required services by the patients.

In order to demonstrate the applicability of the suggested model, a small example was first run to demonstrate the validity of the model. In addition, the interpretation of the results showed the importance of the proposed model. Several sets of problems including large, medium and small-sized instances were used to confirm the efficiency of the new model for problems of various sizes. Moreover, to analyze the effect of important parameters of the problem, a sensitivity analysis was executed on the parameter of the required features for the services. Finally, to help the decision makers to handle their limited resources more effectively, some managerial insights were presented.

As a direction for future research, considering the problem in stochastic situations could be an attractive idea. Employing exact techniques like the branch and cut (B&C) method could be another direction to solve the suggested model. Moreover, when exact approaches cannot be applied to large-sized problems, using meta-heuristic and heuristic approaches could be worthwhile. The home health care problem with considering various existing stakeholders' proposes could be viewed from various perspectives. Thus, multi-objective optimization approaches can be applied to cover different objectives simultaneously. Finally, developing the model using time-dependent traveling times and incorporating nursing education times in the planning could be other streams.

## References

Akjiratikarl, C.; Yenradee, P.; & Drake, P. R. (2007). PSO-Based Algorithm for Home Care Worker Scheduling in the UK. *Computers & Industrial Engineering,* 53(4), 559-583 .

Bazirha, M.; Kadrani, A.; & Benmansour, R. (2020) Scheduling Optimization of the Home Healh Care Problem with Stochastic Travel and Care Times, *Gol International Conference on Logistics Operations Management*, DOI: 10.1109/GOL49479.2020.9314717.

Bazirha, M.; Kadrani, A.; & Benmansour, R. (2021). Stochastic Home Health Care Routing and Scheduling Problem with Multiple Synchronized Services. *Annals of Operations Research*, published online, DOI: 10.1007/s10479-021-04222-w.

Bazirha, S. (2022) Optimization of the Stochastic Home Health Care Routing and Scheduling Problem with Multiple Hard Time Windows. International Journal of Supply and Operatios Management, published online; DOI: 10.22034/IJSOM.2021.109079

Bertels, S.; & Fahle, T. (2006). A Hybrid Ssetup for a Hybrid Scenario: Ccmbining Heuristics for the Home Health Care Problem. *Computers & Operations Research, 33*(10), 2866-2890 .

Carello, G.; Lanzarone, E.; & Mattia, S. (2018). Trade-off between Stakeholders' Goals in the Home Care Nurse-to-Patient Assignment Problem. *Operations Research for Health Care,* 16, 29-40 .

Decerle, J.; Grunder, O.; El Hassani, A. H.; & Barakat, O. (2018). A Memetic Algorithm for a Home Health Care Routing and Scheduling Problem. *Operations Research for Health Care,* 16, 59-71 .

Demirbilek, M.; Branke, J.; & Strauss, A. (2019). Dynamically Accepting and Scheduling Patients for Home Healthcare. *Health Care Management Science,* 22(1), 140-155 .

Di Mascolo, M.; Martinez, C.; & Espinouse, M.-L. (2021). Routing and Scheduling in Home Health Care. *Computers & Industrial Engineering*, 158, 1 – 48.

Euchi, J.; Zidi, S.; & Laouamer, L. (2021). A New Distributed Optimization Approach for Home Healthcare Routing and Scheduling Problem, *Decision Science Letters*, 10 (3), 2021, 217 – 230.

Eveborn, P.; Flisberg, P.; & Rönnqvist, M. (2006). Laps Care—an Operational System for Staff Planning of Home Care. *European Journal of Operational Research,* 171(3), 962-976 .

Fathollahi-Fard, A. M.; Hajiaghaei-Keshteli, M.; & Tavakkoli-Moghaddam, R. (2018). A Bi-objective Green Home Health Care Routing Problem. *Journal of Cleaner Production,* 200, 423-443 .

Fernandez, A.; Gregory, G., Hindle, A., & Lee, A. (1974). A Model for Community Nursing in a Rural County. *Journal of the Operational Research Society,* 25(2), 231-239 .

Fikar, C.; & Hirsch, P. (2018). Evaluation of Trip and Car Sharing Concepts for Home Health Care Services. *Flexible Services and Manufacturing Journal,* 30(1-2), 78-97 .

Ghiasi, F.G.; Yazdani, M.; Vahdani, B.; & Kazemi, A. (2021). Multi-Depot Home Health Care Routing and Scheduling Problem with Multimodal Transportation: Mathematical Model and Solution Methods, *Scientia Iranica*, published online; DOI: 10.24200/SCI.2021.57338.5183

Grenouilleau, F.; Legrain, A.; Lahrichi, N., & Rousseau, M. (2019) A Set Partioning Heuristic for the Home Health Care Routing and Scheduling Problem, *European Journal of Operational Research*, 275 (1), 295 - 303.

Hindle, T.; Hindle, A.; & Spollen, M. (2000). Resource Allocation Modelling for Home-Based Health and Social Care Services in Areas Having Differential Population Density Levels: a Case Study in Northern Ireland. *Health Services Management Research,* 13(3), 164-169 .

Hindle, T.; Hindle, G.; & Spollen, M. (2009). Travel-Related Costs of Population Dispersion in the Provision of Domiciliary Care to the Elderly: a Case Study in English Local Authorities. *Health Services Management Research,* 22(1), 27-32 .

Issabakhsh, M.; Hosseini-Motlagh, S.-M.; Pishvaee, M.-S.; & Saghafi Nia, M. (2018). A Vehicle Routing Problem for Modeling Home Healthcare: a Case Study. *International Journal of Transportation Engineering,* 5(3), 211-228 .

Kohdabandeh, P.; Rafiee, M.; Kayvanfar, V.; & Werner, F. (2021). A Bi-objective Home Health Care Routing and Scheduling Model with Considering Nurse Downgrading Costs, *International Journal of Enviromental Research and Public Health*, 18 (3), 900, 24 pages.

Khodaparasti, S.; Bruni, M. E.; Beraldi, P.; Maleki, H.; & Jahedi, S. (2018). A Multi-Period Location-Allocation Model for Nursing Home Network Planning under Uncertainty. *Operations Research for Health Care,* 18, 4-15 .

Lanzarone, E.; & Matta, A. (2014). Robust Nurse-to-Patient Assignment in Home Care Services to Minimize Overtimes under Continuity of Care. *Operations Research for Health Care,* 3(2), 48-58 .

Lin, C.-C.; Hung, L.-P.; Liu, W.-Y.; & Tsai, M.-C. (2018). Jointly Rostering, Routing, and Rerostering for Home Health Care Services: A Harmony Search Approach with Genetic, Saturation, Inheritance, and Immigrant Schemes. *Computers & Industrial Engineering,* 115, 151-166 .

Liu, M.; Yang, D.; Su, Q.; & Xu, L. (2018). Bi-objective Approaches for Home Healthcare Medical Team Planning and Scheduling Problem. *Computational and Applied Mathematics,* 37(4), 4443-4474 .

Liu, R.; Xie, X.; & Garaix, T. (2014). Hybridization of Tabu Search with Feasible and Infeasible Local Searches for Periodic Home Health Care Logistics. *Omega,* 47, 17-32 .

Liu, R.; Yuan, B.; & Jiang, Z. (2018). A Branch-and-Price Algorithm for the Home-Caregiver Scheduling and Routing Problem with Stochastic Travel and Service Times. *Flexible Services and Manufacturing Journal*, 1-23 .

Mankowska, D. S.; Meisel, F.; & Bierwirth, C. (2014). The Home Health Care Routing and Scheduling Problem with Interdependent Services. *Health Care Management Science,* 17(1), 15-30 .

Nasir, J.; & Dang, C. (2018). Solving a More Flexible Home Health Care Scheduling and Routing Problem with Joint Patient and Nursing Staff Selection. *Sustainability,* 10(1), 148 .

Nasir, J. A., Hussain, S., & Dang, C. (2018). An Integrated Planning Approach Towards Home Health Care, Telehealth and Patients Group Based Care. *Journal of Network and Computer Applications,* 117, 30-41 .

Rodriguez, C., Garaix, T., Xie, X., & Augusto, V. (2015). Staff dimensioning in homecare services with uncertain demands. *International Journal of Production Research,* 53(24), 7390-7410.

Shanejat-Bushehri, S.; Tavakkoli-Moghaddam R.; Momen, S.; Gasemkhani, A.; Tavakkoli-Moghaddam, H. (2019). Home Health Care Routing and Scheduling Problem Considering Temporal Dependencies and Perishability with Simultaneous Pickup and Delivery, *IFAC Papers Online*, 52, No. 13, 118 – 123.

Shi, Y.; Boudouh, T.; & Grunder, O. (2017). A hybrid genetic algorithm for a home health care routing problem with time window and fuzzy demand. *Expert Systems with Applications,* 72, 160-176 .

Shi, Y.; Boudouh, T.; Grunder, O.; & Wang, D .(2018) .Modeling and Solving Simultaneous Delivery and Pick-up Problem with Stochastic Travel and Service Times in Home Health Care. *Expert Systems with Applications,* 102, 218-233 .

Xiang, T.; Li, Y; & Szeto, W.Y. (2021). The Daily Routing and Scheduling of Home Health Care: Based on Costs and Participant's Preference Satisfaction. *International Transactions in Operational Research*, Published online, DOI: 10.1111/itor.1304.

Yuan, B.; Liu, R.; & Jiang, Z. (2015). A Branch-and-Price Algorithm for the Home Health Care Scheduling and Routing Problem with Stochastic Service Times and Skill Requirements. *International Journal of Production Research,* 53(24), 7450-7464 .

Zhang, T.; Yang, X.; Chen, Q.; Bai, L.; & Chen, W. (2018). *Modified ACO for Home Health Care Scheduling and Routing Problem in Chinese Communities.* Paper presented at the 2018 IEEE 15th International Conference on Networking, Sensing and Control (ICNSC).